# DISCUSSION OF "ANALYSIS OF VARIANCE—WHY IT IS MORE IMPORTANT THAN EVER" BY A. GELMAN

By Tue Tjur

*Copenhagen Business School*

I am not sure that I agree with the statement that ANOVA is more important than ever. On the contrary, I think that the development during the last 40–50 years has somewhat scaled down the importance of this topic, by separating the computational aspect from the model aspect. Many of the concepts in classical ANOVA for balanced designs are related to the computations. Of course, it is still easier to do the computations in the balanced case, and balancedness also implies other advantages such as maximal efficiency and exact distributions instead of approximations for the mixed models. But the availability of methods for handling of linear models and mixed models in unbalanced designs has changed the focus. Today, I believe, we are more inclined to think of these models as examples (though very important examples) of statistical models, whereas in the classical approach one could hardly mention ANOVA and multiple regression in the same course or textbook.

To me, modern statistics [as opposed, e.g., to the approach taken by Cochran and Cox (1957)] is characterized by the ultimate focusing on the statistical model as the central object. And this brings me to the main topic of my comment to this article, which is the theory or method presented in Section 3. I must admit that I am rather confused here and that I have not been able to understand much of what is going on. The reason for this is, as I see it, that it is not clear at all what the statistical model is. The basic idea seems to be to let all effects enter formally as random effects. But since the method is claimed to be able to handle fixed effects as well (and even to make the comparisons automatically with the correct standard deviations), there must be something I have missed. The random model is only a tool, it is obviously not the model we want to analyze.

Intuitively, it is not difficult to see that there is some element of truth in this approach. For example, in the machine-treatment example, where treatments are confounded with machines, it is certainly correct that the interesting part of the analysis is equivalent to a simple one-way analysis of the 20 machine averages. But what is not at all clear to me is what the method actually does (a detailed operational description), when it works and why it works. I can see no way of proving this without an explicit statement of the model(s) that we actually want to analyze. My guess is that the validity of this method (whatever it is) can only be proved under assumptions about balancedness and orthogonality. Even here there may be problems, since it is not obvious how a phenomenon like partial confounding of a treatment effect with a block effect can be handled. Probably by the introduction of pseudo factors, but where do they come in?

This is all rather negative, and I would have liked to be more positive because I think one of the declared purposes, to make split-plot and other analyses more understandable to students, is an important one. However, my experience here is that the best way of making these things understandable is to focus on the model rather than the design. The analysis of a split-plot design should, in my opinion, be regarded as no more and no less than the analysis of a mixed model. The implications of balancedness (considerable simplification of the computations, exact confidence intervals for contrasts, exact distributions of test statistics, etc.) are important, but irrelevant to the understanding of the statistical model and the interpretation of its parameters.

The Statistics Group
Copenhagen Business School
Solbjerg Plads 3
DK-2000 Frederiksberg
Denmark
e-mail: tt.mes@cbs.dk